\documentclass[12pt,twosided,reqno]{amsart}
\usepackage{amsmath}
\usepackage{amsthm}
\usepackage{amsfonts}
\usepackage{amssymb}
\usepackage{color}

\theoremstyle{theorem}

\theoremstyle{definition}

\begin{document}

\vskip .5cm
\title{Some characterizations of spheres and elliptic paraboloids II}
 \vskip0.5cm
 \author{Dong-Soo Kim$^1$ and Young Ho Kim$^2$}
\address{\newline Department of Mathematics, Chonnam National University, Kwangju
500-757, Korea \\ \newline
Department of Mathematics, Kyungpook
National University, Taegu 702-701, Korea} \email{
dosokim@chonnam.ac.kr and yhkim@knu.ac.kr}
\date{}
\dedicatory{}
\thanks{
\newline\indent
     \newline\indent { $^1$  supported by Basic Science Research Program through the National Research Foundation of Korea (NRF)
     funded by the Ministry of Education, Science and Technology (2010-0022926).}
\newline\indent { $^2$   supported by Basic Science Research Program through the National Research Foundation of Korea (NRF)
 funded by the Ministry of Education, Science and Technology (2010-0007184).}}

\begin{abstract}
We show some characterizations of hyperspheres  in the $(n+1)$-dimensional Euclidean space ${\Bbb E}^{n+1}$ with intrinsic and extrinsic properties such as the $n$-dimensional area of the sections cut off by hyperplanes, the $(n+1)$-dimensional volume of regions between parallel hyperplanes, and the $n$-dimensional surface area of regions between parallel hyperplanes.
We also establish two  characterizations of elliptic paraboloids  in the $(n+1)$-dimensional Euclidean space
 ${\Bbb E}^{n+1}$ with the $n$-dimensional area of the sections cut off by hyperplanes and the $(n+1)$-dimensional volume of regions between parallel hyperplanes.
For further study, we suggest a few open problems.

\end{abstract}

\vskip 1cm

\date{}
\maketitle
\section {Introduction} \vskip 0.3cm

Let  $S^n(a)$ be a hypersphere with radius $a$ in the Euclidean space
 ${\Bbb E}^{n+1}$.  For a fixed point $p \in S^n(a)$ and  for a sufficiently small $t>0$,
 let's denote by $\Phi$ a hyperplane  parallel to
 the tangent space $\Psi$ of $S^n(a)$ at $p$ with distance $t$
 which intersects $S^n(a)$.

  We denote by $A_p(t), V_p(t)$ and $S_p(t)$
 the $n$-dimensional area of the section in $\Phi$ enclosed by $\Phi \cap S^n (a)$,
 the $(n+1)$-dimensional volume of the region bounded by the sphere and the plane $\Phi$ and the $n$-dimensional surface area of the region of $S^n(a)$  between the two planes $\Phi$ and $\Psi$,
respectively.

\vskip 0.3cm
 Then, for a sufficiently small $t>0$, we can have the following properties of the sphere $S^n(a)$.
\vskip 0.3cm

 \noindent $(A)$: The $n$-dimensional area $A_p(t)$ of the section  is independent of the point $p$.

 \noindent $(V)$: The $n$-dimensional volume  $V_p(t)$ of the region   is independent of the point $p$.

  \noindent $(S)$: The $n$-dimensional surface area  $S_p(t)$ of the region   is independent of the point $p$.

\vskip 0.3cm

If $n=2$, Archimedes proved that $S_p(t)=2\pi at$ holds for   $S^2(a)$ (\cite{se}, p.78).
For a differential geometric proof, see Archimedes' Theorem (\cite{pr}, pp.116-118).
\smallskip

 Conversely, it is natural to ask the following question:
\vskip 0.30cm

 \noindent {\bf Question 1.} ``Are there any other hypersurfaces in Euclidean space which  satisfy the above properties?"
\vskip 0.30cm

For  the case of $n=2$ about the property $(S)$, the authors answered negatively  as follows   (\cite{kk1}) (See also \cite{b} and \cite{s}.):
\vskip 0.30cm

 \noindent {\bf Proposition 2.}
 Let $M$ be a closed and convex surface in the 3-dimensional Euclidean space
 ${\Bbb E}^{3}$. Suppose that  $M$ satisfies the condition:

 \noindent  (S) $S_p(t)=\phi(t)$,
 which depends only on $t$.

Then $M$ is a round sphere.
 \vskip0.3cm

In this article, first, we study convex hypersurfaces $M$ in the $(n+1)$-dimensional Euclidean space
 ${\Bbb E}^{n+1}$ which satisfy the above mentioned properties. For a point $p\in M \subset {\Bbb E}^{n+1}$,
$A_p(t), V_p(t)$ and $S_p(t)$ are defined as above.

In Section 3, as a result, we prove the following:
\vskip 0.30cm

 \noindent {\bf Theorem 3.}
 Let $M$ be a complete  and convex hypersurface in the $(n+1)$-dimensional Euclidean space
 ${\Bbb E}^{n+1}$. Suppose that $M$ satisfies one of the following conditions.

\vskip 0.3cm

 \noindent $(A)$: The $n$-dimensional area $A_p(t)$ of the section  is independent of the point $p\in M$.

 \noindent $(V)$: The $(n+1)$-dimensional volume  $V_p(t)$ of the region   is independent of the point $p\in M$.

  \noindent $(S)$: The $n$-dimensional surface area  $S_p(t)$ of the region   is independent of the point $p\in M$.

\vskip 0.3cm
Then the hypersurface $M$ is a round hypersphere $S^n(a)$.

 \vskip 0.3cm

Second, suppose that $M$ is a smooth convex  hypersurface  in the $(n+1)$-dimensional Euclidean space
 ${\Bbb E}^{n+1}$ defined  by the graph of a convex function  $f:{\Bbb R}^{n}\rightarrow {\Bbb R}$.
 For a fixed point $p=(x,f(x)) \in M$ and  for a real number $k>0$,
 consider a hyperplane $\Phi$ through $v=(x, f(x)+k)$ which is parallel to
 the tangent hyperplane $\Psi$ of $M$ at $p$.

  We denote by $A_p^*(k), V_p^*(k)$ and $S_p^*(k)$
 the $n$-dimensional area of the section in $\Phi$ enclosed by $\Phi \cap M$,
 the $(n+1)$-dimensional volume of the region bounded by $M$ and the hyperplane $\Phi$, and
the $n$-dimensional surface area of the region of $M$  between the two hyperplanes $\Phi$ and $\Psi$,
respectively.

\vskip0.3cm
For elliptic paraboloids,  in \cite{kk1} the authors proved the following.
\vskip 0.50cm

 \noindent {\bf Proposition 4.}
  Let $M$ be a smooth convex  hypersurface  in the $(n+1)$-dimensional Euclidean space
 ${\Bbb E}^{n+1}$ defined  by the graph of a convex function  $f:{\Bbb R}^{n}\rightarrow {\Bbb R}$.
  Then $M$ is an elliptic paraboloid if and only if   it satisfies  the following condition:

  \noindent $(L)$: The $(n+1)$-dimensional volume $V_p^*(k)$ is $ak^{(n+2)/2}$ for some constant $a$
  which depends only on the hypersurface $M$.

 \vskip0.3cm

In Section 3, we generalize  Proposition 4 as follows.

 \vskip 0.3cm

 \noindent {\bf Theorem 5.}
 Let $M$ be a smooth convex  hypersurface  in the $(n+1)$-dimensional Euclidean space
 ${\Bbb E}^{n+1}$ defined  by the graph of a convex function  $f:{\Bbb R}^{n}\rightarrow {\Bbb R}$.
  Then $M$ is an elliptic paraboloid if and only if   it satisfies one of the following conditions:

 \vskip 0.3cm
 \noindent $(V^*)$:  $V_p^*(k)$    is  a nonnegative function $\phi (k)$,
 which depends only on $k$.

 \noindent $(A^*)$:
  $A_p^*(k)/W(p)$ is  a nonnegative function $\psi (k)$,
 which depends only on $k$.

 Here, we denote $p=(x,f(x))$ and  $W(p)=\sqrt{1+|\nabla f(x)|^2}$,  where $\nabla f$ is the
  gradient of $f$.

\vskip 0.3cm
In view of the conditions in Theorem 5,  it is reasonable to ask  the following question.

\vskip 0.30cm

 \noindent {\bf Question 6.} Which hypersurfaces  satisfy the following condition $(S^*)$?

  \noindent $(S^*)$: $S_p^*(k)/W(p)$ is  a nonnegative function $\eta (k)$, which depends only on $k$.
  \vskip 0.30cm

Finally, using harmonic function theory (\cite{a}), we answer Question 6 negatively as follows.
\vskip 0.3cm

\noindent {\bf Theorem 7.} Let $M$ be a smooth convex  hypersurface  in the $(n+1)$-dimensional Euclidean space
 ${\Bbb E}^{n+1}$ defined  by the graph of a convex function  $f:{\Bbb R}^{n}\rightarrow {\Bbb R}$.
Then $M$ does not satisfy condition $(S^*)$.

 \vskip0.3cm

In this paper, in order to prove  our theorems,
 we prove a lemma (Lemma 8), extending a lemma in \cite{kk1}, about a new meaning of Gauss-Kronecker curvature $K(p)$
of convex hypersurface $M$ at a point $p\in M$ in three ways.

  \vskip 0.3cm

  We now state some questions for further study  as follows.

    \vskip 0.3cm

\noindent {\bf  Question A.}
 Let $M$ be a  convex (not complete)  hypersurface in the $(n+1)$-dimensional Euclidean space
 ${\Bbb E}^{n+1}$. Suppose that $M$ satisfies one of the  conditions in Theorem 3.
Then, is it an open part of a round hypersphere $S^n(a)$?
\vskip 0.3cm

An elliptic paraboloid satisfies the following conditions. For a proof, see the proof of Theorem 5, which is given in  Section 3.

\noindent $(V^{**})$: $V_p(t)=C(p)t^{(n+2)/2}$, where $C(p)$ is a function of $p\in M$.

\noindent $(A^{**})$: $A_p(t)=D(p)t^{n/2}$, where $D(p)$ is a function of $p\in M$.

Due to (2.5), the above two conditions are equivalent.

\vskip 0.3cm

 \noindent {\bf Question B.}
 Let $M$ be a complete  and convex hypersurface in the $(n+1)$-dimensional Euclidean space
 ${\Bbb E}^{n+1}$, which is not necessarily a graph of a function.
  Suppose that $M$ satisfies the condition $(V^{**})$.
 Then,  is it  an elliptic paraboloid?

\vskip 0.3cm
For $n=1$,  Question A is true because the plane curvature is a nonzero constant.
In [6], the authors  answered Question B for $n=1$, affirmatively.
\vskip 0.3cm

Throughout this article, all objects are smooth and connected, otherwise mentioned.

\section {Preliminaries}
\vskip 0.3cm
Suppose that $M$ is a smooth convex  hypersurface  in the $(n+1)$-dimensional Euclidean space
 ${\Bbb E}^{n+1}$.
 For a fixed point $p \in M$ and  for a sufficiently small $t>0$, consider a hyperplane $\Phi$ parallel to
 the tangent hyperplane $\Psi$ of $M$ at $p$ with distance $t$
 which intersects $M$.

 We denote by $A_p(t), V_p(t)$ and $S_p(t)$
 the $n$-dimensional area of the section in $\Phi$ enclosed by $\Phi \cap M$,
 the $(n+1)$-dimensional volume of the region bounded by the hypersurface and the hyperplane $\Phi$ and
the $n$-dimensional surface area of the region of $M$  between the two hyperplanes $\Phi$ and $\Psi$,
respectively.

Now, we may introduce a coordinate system $(x,z)=(x_1,x_2, \cdots , x_n, z)$
 of  ${\Bbb E}^{n+1}$ with the origin $p$, the tangent  space of $M$ at $p$ is the hyperplane $z=0$.
 Furthermore, we may assume that $M$ is locally  the graph of a non-negative convex  function $f:{\Bbb R}^n\rightarrow {\Bbb R}$.

Then, for a sufficiently small $t>0$ we have
\smallskip
   \begin{equation}\tag{2.1}
   \begin{aligned}
   A_p(t)&= \iint _{f(x)<t}1dx,
    \end{aligned}
   \end{equation}
   \begin{equation}\tag{2.2}
   \begin{aligned}
   V_p(t)&=\iint _{f(x)<t}\{t-f(x)\}dx
    \end{aligned}
   \end{equation}
   and
 \begin{equation}\tag{2.3}
   S_p(t) = \iint _{f(x)<t}\sqrt{1+|\nabla f|^2}dx,
   \end{equation}
 where $x=(x_1,x_2, \cdots , x_n)$, $dx=dx_1dx_2 \cdots dx_n$ and $\nabla f$ denotes the gradient vector of the function $f$.

Note that  we also have
 \begin{equation}\tag{2.4}
   \begin{aligned}
   V_p(t)&=\iint _{f(x)<t}\{t-f(x)\}dx\\
   &=\int _{z=0}^{t}\{\iint_{f(x)<z}1dx\}dz.
       \end{aligned}
   \end{equation}
Hence, together with the fundamental theorem of calculus,  (2.4) shows that
\begin{equation}\tag{2.5}
   \begin{aligned}
  V_p'(t)=\iint _{f(x)<t}1dx=A_p(t).
    \end{aligned}
   \end{equation}

First of all, we  prove
\vskip 0.50cm

 \noindent {\bf Lemma 8.} Suppose that the Gauss-Kronecker curvature
 $K(p)$ of $M$ at $p$ is positive with respect to the upward unit normal to $M$.
Then we have the following:

\noindent 1)
 \begin{equation}\tag{2.6}
 \lim_{t \rightarrow 0}\frac{1}{(\sqrt{t})^{n}}A_p(t)= \frac{(\sqrt{2})^{n}\omega_n}{\sqrt{K(p)}},
 \end{equation}

\noindent 2)
 \begin{equation}\tag{2.7}
 \lim_{t \rightarrow 0}\frac{1}{(\sqrt{t})^{n+2}}V_p(t)= \frac{(\sqrt{2})^{n+2}\omega_n}{(n+2)\sqrt{K(p)}},
 \end{equation}

\noindent 3)
 \begin{equation}\tag{2.8}
 \lim_{t \rightarrow 0}\frac{1}{(\sqrt{t})^{n}}S_p(t)= \frac{(\sqrt{2})^{n}\omega_n}{\sqrt{K(p)}},
 \end{equation}
where $\omega_n$ denotes the volume of  the $n$-dimensional unit ball.
\smallskip

\vskip0.3cm
\noindent {\bf Proof.} We denote by $x$ the column vector $(x_1, x_2, \cdots , x_n)^t$. Then
for a symmetric $n\times n$ matrix $A$, we have from the Taylor's  formula of $f(x)$ as follows:
 \begin{equation}\tag{2.9}
 f(x)= x^tAx + f_3(x),
  \end{equation}
where $f_3(x)$ is an $O(|x|^3)$  function.
Then the Hessian matrix of $f$ at the origin is given by
$D^2f(0)=2A.$
Hence, for the upward unit normal to $M$ we have
 \begin{equation}\tag{2.10}
 K(p)=\det D^2f(0) = 2^n \det A.
  \end{equation}
By the assumption, we see that every eigenvalue of $A$ is positive.
Thus,  there exists a nonsingular symmetric matrix $B$ satisfying
 \begin{equation}\tag{2.11}
 A=B^tB,
 \end{equation}
 where $B^t$ denotes the transpose of $B$.
Therefore, we obtain
 \begin{equation}\tag{2.12}
 f(x)= |Bx|^2+f_3(x).
  \end{equation}

We consider the decomposition of $S_p(t)$ as follows:
 \begin{equation}\tag{2.13}
 \begin{aligned}
 S_p(t)&= A_p(t) + N_p(t),
   \end{aligned}
  \end{equation}
  where
 \begin{equation}\tag{2.1}
 \begin{aligned}
 A_p(t)&=\iint_{f(x)<t}1dx
  \end{aligned}
  \end{equation}
and
 \begin{equation}\tag{2.14}
 \begin{aligned}
  N_p(t)&=\iint_{f(x)<t}(\sqrt{1+|\nabla f|^2}-1)dx.
   \end{aligned}
  \end{equation}

First, we show (2.6) as follows.
We let $t=\epsilon^2$ and $x=\epsilon y$. Then (2.1) gives
 \begin{equation} \tag{2.15}
  \begin{aligned}
  \frac{1}{(\sqrt{t})^{n}}A_p(t)= \frac{1}{(\sqrt{t})^{n}}\iint_{f(x)<t}1dx =\iint_{|By|^2+\epsilon g_3(y)<1}1dy,
 \end{aligned}
  \end{equation}
where $g_3(y)$ is an $O(|y|^3)$  function. As $\epsilon \rightarrow 0$,
it follows from (2.15) that
\begin{equation} \tag{2.16}
  \begin{aligned}
 \lim_{t\rightarrow 0}  \frac{1}{(\sqrt{t})^{n}}A_p(t)=\iint_{|By|^2<1}1dy.
 \end{aligned}
  \end{equation}

If we let $w=By$, then from (2.16) we get
\begin{equation} \tag{2.17}
  \begin{aligned}
 \lim_{t\rightarrow 0}  \frac{1}{(\sqrt{t})^{n}}A_p(t)
=\frac{1}{|\det B|}\iint_{|w|<1}1dw=\frac{\omega_n}{|\det B|}.
 \end{aligned}
 \end{equation}
 Hence, it follows from (2.10) and (2.11) that (2.6) holds.

Together with (2.5) and (2.6), L'Hospital's rule implies (2.7).

In order to prove (2.8), it suffices to show that
\begin{equation} \tag{2.18}
  \begin{aligned}
 \lim_{t\rightarrow 0}  \frac{1}{(\sqrt{t})^{n}}N_p(t)=0.
 \end{aligned}
 \end{equation}
 Note that the following inequality holds
 \begin{equation} \tag{2.19}
  \begin{aligned}
N_p(t)\le \frac{1}{2}\iint_{f(x)<t}|\nabla f|^2dx.
 \end{aligned}
 \end{equation}
The function $f$ satisfies
 \begin{equation} \tag{2.20}
  \begin{aligned}
|\nabla f(x)|^2=4|Ax|^2+h_2(x),
 \end{aligned}
 \end{equation}
where $h_2(x)$ is an $O(|x|^2)$ function.
Thus, there exists a  positive constant $C$ satisfying  in a neighborhood of the origin
 \begin{equation} \tag{2.21}
  \begin{aligned}
|\nabla f(x)|^2 \le C|x|^2.
 \end{aligned}
 \end{equation}

In the same argument as above, putting  $t=\epsilon^2$ and $x=\epsilon y$,
it follows from (2.19) and (2.21) that
 \begin{equation} \tag{2.22}
  \begin{aligned}
  \frac{1}{(\sqrt{t})^{n}}N_p(t)\le  \frac{C\epsilon^2}{2}\iint_{|By|^2+\epsilon g_3(y)<1}|y|^2 dy.
 \end{aligned}
  \end{equation}
Since the integral of the right side  in (2.22) tends to a constant as $\epsilon \rightarrow 0$,
 by letting $t \rightarrow 0$ in (2.22), we get (2.18). Together with (2.6) and (2.13), (2.18) shows that (2.7) holds.
 This completes the proof. $\qquad$ $\square$
 \vskip 0.50cm

\section {Proofs of theorems}
\vskip 0.3cm

In this section, first, we prove Theorem 3.

Let $M$ be a complete  and convex hypersurface in the $(n+1)$-dimensional Euclidean space
 ${\Bbb E}^{n+1}$. Then
 the Gauss-Kronecker curvature $K(p)$ of $M$ with respect to the inward unit normal to $M$
 is nonnegative and positive somewhere.

  Suppose that $M$ satisfies one of the conditions in Theorem 3.
 Then Lemma 8 shows that the Gauss-Kronecker curvature $K(p)$ is  constant on the nonempty open set
 $\Omega=\{p\in M|K(p)>0\}$. Hence the continuity of $K$ implies that $\Omega=M$, that is, $K(p)$ is  constant on $M$.
 Thus,  it  follows from  Theorem 7.1 in [9] or the main theorem in [3]
 that $M$ is a round hypersphere.

 The converse is obvious.

 \vskip 0.50cm
Second,  we give a proof of Theorem 5.

Suppose that $M$ is a smooth convex  hypersurface  in the $(n+1)$-dimensional Euclidean space
 ${\Bbb E}^{n+1}$ defined  by the graph of a convex function  $f:{\Bbb R}^{n}\rightarrow {\Bbb R}$.
We fix a  point
$p=(x,f(x))$ in $M$. For a positive constant $k$, consider a hyperplane  $\Phi$ through $v=(x,f(x)+k),$
which is parallel to the tangent hyperplane $\Psi$ to $M$ at $p$.
Let $W(x)=\sqrt{1+|\nabla f(x)|^2}$.
Then  for a constant $t$ with $k=tW$, we have  $V_p^*(k)=V_p(t)$, $A_p^*(k)=A_p(t)$ and $S_p^*(k)=S_p(t)$.

\vskip 0.30cm
\noindent 1) The condition $(V^*)$  shows that $V_p(t)=\phi(k)$, which is independent of $p\in M$.
Hence we have
\begin{equation} \tag{3.1}
  \begin{aligned}
\frac{V_p(t)}{(\sqrt{t})^{n+2}}= \frac{\phi(k)}{(\sqrt{k})^{n+2}}(\sqrt{W})^{n+2}.
 \end{aligned}
 \end{equation}

Thus,  Lemma 8 implies that
\begin{equation}\tag{3.2}
   \begin{aligned}
  \lim_{k\rightarrow 0}\frac{\phi(k)}{(\sqrt{k})^{n+2}}=
   W^{-(n+2)/2}\lim_{t\rightarrow 0} \frac{V_p(t)}{(\sqrt{t})^{n+2}}= \frac{(\sqrt{2})^{n+2}\omega_n}{(n+2)\sqrt{K(p)}}  W^{-(n+2)/2}.
    \end{aligned}
   \end{equation}
If we denote by $\alpha$ the limit of the left side of (3.2), which is independent of $p$,  then we have
\begin{equation}\tag{3.3}
   \begin{aligned}K(p)=\frac{2^{n+2}\omega_n^2}{\alpha^2(n+2)^2W(x)^{n+2}}.
    \end{aligned}
   \end{equation}
Since the Gauss-Kronecker curvature $K(p)$ of $M$ at $p$ is given by ([12], p.93)
\begin{equation}\tag{3.4}
   \begin{aligned}
   K(p)=\frac{\det D^2f(x)}{W^{n+2}},
    \end{aligned}
   \end{equation}
it follows from (3.3) that the determinant $\det D^2f(x)$
of the Hessian of $f(x)$ is a positive constant.
The continuity of $\det D^2f(x)$ shows that it is a positive constant on the whole space ${\Bbb R}^n$.
Thus $f(x)$ is a globally defined quadratic polynomial (\cite{j,po}), and hence $M$ is an elliptic paraboloid.

\vskip 0.30cm
\noindent 2) The condition $(A^*)$ shows that $A_p^*(k)/W(p)=\psi (k)$, which is independent of $p\in M$.
Hence we have
\begin{equation} \tag{3.5}
  \begin{aligned}
\frac{A_p(t)}{(\sqrt{t})^{n}}= \frac{\psi(k)}{(\sqrt{k})^{n}}(\sqrt{W})^{n+2}.
 \end{aligned}
 \end{equation}

Thus,  Lemma 8 implies that
\begin{equation}\tag{3.6}
   \begin{aligned}
  \lim_{k\rightarrow 0}\frac{\psi(k)}{(\sqrt{k})^{n}}=
   W^{-(n+2)/2}\lim_{t\rightarrow 0} \frac{A_p(t)}{(\sqrt{t})^{n}}= \frac{(\sqrt{2})^{n}\omega_n}{\sqrt{K(p)}}  W^{-(n+2)/2}.
    \end{aligned}
   \end{equation}
If we denote by $\beta$ the limit of the left side of (3.6), which is independent of $p$,  then we have
\begin{equation}\tag{3.7}
   \begin{aligned}K(p)=\frac{2^{n}\omega_n^2}{\beta^2W(x)^{n+2}}.
    \end{aligned}
   \end{equation}
Hence, as in the proof of Case 1), we see that  $M$ is an elliptic paraboloid.

\vskip 0.30cm

 This completes the proof of the if  part of Theorem 5.

\vskip 0.30cm
Conversely, consider an elliptic paraboloid  $M:z=f(x)=\Sigma_{i=1}^na_i^2x_i^2, a_i>0$,
 a tangent hyperplane $\Psi$ to $M$ at a fixed point $p=(x,z) \in M$,
 a hyperplane $\Phi$ through $v=(x, z+k), k>0$ which is parallel to
 the tangent hyperplane $\Psi$ of $M$ at $p$. Then the proof of Theorem 5 of [5] shows that
 \begin{equation}\tag{3.8}
   \begin{aligned}
   V_p^*(k)=\alpha_nk^{(n+2)/2},\quad \alpha_n=\frac{2\sigma_{n-1}}{n(n+2)a_1a_2\cdots a_n},
    \end{aligned}
   \end{equation}
where $\sigma_{n-1}$ denotes the surface area of the $(n-1)$-dimensional unit sphere.
Since $V_p^*(k)=V_p(t)$ with $k=tW$, we get
 \begin{equation}\tag{3.9}
   \begin{aligned}
   V_p(t)=\alpha_nW(p)^{(n+2)/2}t^{(n+2)/2}.
    \end{aligned}
   \end{equation}

Hence, it follows from (2.5) that
  \begin{equation}\tag{3.10}
   \begin{aligned}
   A_p(t)=\frac{n+2}{2}\alpha_nW(p)^{(n+2)/2}t^{n/2}
    \end{aligned}
   \end{equation}
 and
  \begin{equation}\tag{3.11}
   \begin{aligned}
   A_p^*(k)=\frac{n+2}{2}\alpha_nW(p)k^{n/2}.
    \end{aligned}
   \end{equation}
 Thus, (3.8) and (3.11) show that an elliptic paraboloid
 $M$ satisfies conditions $(V^*)$ and $(A^*)$ in Theorem 5, respectively.

 This completes the proof of the only if  part of Theorem 5.

\vskip 0.30cm
Finally, we prove Theorem 7.

Suppose that a smooth convex  hypersurface  $M$ in the $(n+1)$-dimensional Euclidean space
 ${\Bbb E}^{n+1}$ defined  by the graph of a convex function  $f:{\Bbb R}^{n}\rightarrow {\Bbb R}$
 satisfies the condition $(S^*)$. Then, as in the proof of Theorem 5,
 we can prove that $M$ is an elliptic paraboloid given by $z=f(x)=\Sigma_{i=1}^na_i^2x_i^2, a_i>0$.

Consider
 a hyperplane $\Phi$ intersecting $M$, a point $p=(p_1, \cdots , p_n, \Sigma_{i=1}^na_i^2p_i^2)\in M$
 where the tangent hyperplane $\Psi$ of $M$ is parallel to $\Phi$,
 and a point $v$ where the line through $p$ parallel to the $z$-axis meets $\Phi$ with $||p-v||=k$.
Then we have the following.
\begin{equation}\tag{3.12}
   \begin{aligned}
   \Phi: z&=2a_1^2p_1x_1+\cdots +2a_n^2p_nx_n-(a_1^2p_1^2+\cdots +a_n^2p_n^2)+k, \\
   S_p^*(k)&=\int_{D_p(k)}W(x)dx,\\
   D_p(k)&: \Sigma_{i=1}^na_i^2(x_i-p_i)^2<k,
    \end{aligned}
   \end{equation}
   where $W(x)=\sqrt{1+\Sigma_{i=1}^n4a_i^4x_i^2}.$

By the linear transformation $y_i=a_ix_i, i=1, 2, \cdots , n$, we obtain
\begin{equation}\tag{3.13}
   \begin{aligned}
    S_p^*(k)&=\frac{1}{a_1a_2\cdots a_n}\int_{B_q(\sqrt{k})}V(y)dy,\\
   B_q(\sqrt{k})&: \Sigma_{i=1}^n(y_i-q_i)^2<k,
    \end{aligned}
   \end{equation}
  where $V(y)=\sqrt{1+4\Sigma_{i=1}^na_i^2y_i^2}$ and $q=(a_1p_1, \cdots, a_np_n)$.

 Since  $M$ satisfies the condition $(S^*)$ with $W(p)=V(q)$, by letting $r=\sqrt{k}$,
  it follows from (3.13) that $V(y)$ satisfies the
  following.
\begin{equation}\tag{3.14}
   \begin{aligned}
    \int_{B_q(r)}V(y)dy=V(q)g(r), \quad q\in R^n, r\ge 0,
    \end{aligned}
   \end{equation}
where $ B_q(r)=\{y||y-q|<r\}$ is the ball  of radius $r$ centered at $q$ and $g(r)$ is a function of $r$.

For a function $g=g(r), r\ge 0$, we denote by $C_g$ the set of all functions  $f:{\Bbb R}^{n}\rightarrow {\Bbb R}$ satisfying (3.14).
Then, it is straightforward to show the following (\cite{a}).
\vskip0.3cm

\noindent {\bf Lemma 9.} The set $C_g$ satisfies the following.

\noindent 1) If $g(r)$ is the volume $\omega_nr^n$ of $B_q(r)$  for sufficiently small $r>0$,
then $C_g$ is the set of all harmonic functions on $R^n$.

\noindent 2) If a positive function in $C_g$ has a local maximum (local minimum, respectively),
then $g(r)\le \omega_nr^n$ ($g(r)\ge \omega_nr^n$, respectively) for sufficiently small $r>0$.

\noindent 3) If $f\in C_g$, then every partial derivative of $f$ also belongs to  $C_g$.

\noindent 4) Every linear combination of functions in  $C_g$ also belongs to  $C_g$.

\vskip0.3cm
In order to complete the proof of Theorem 7, we use Lemma 9 as follows. By differentiating, we have
\begin{equation}\tag{3.15}
   \begin{aligned}
  V_{ii}=4a_i^2\frac{(V^2-4a_i^2x_i^2)}{V^3}, \quad i=1,2,\cdots , n,
    \end{aligned}
   \end{equation}
where $V_i$ means the $i$-th partial derivative of $V$,  etc.. Hence, we get
\begin{equation}\tag{3.16}
   \begin{aligned}
  U=\frac{1}{4}\Sigma_{i=1}^{n} (\frac{V_{ii}}{a_{i}^2})=\frac{(n-1)V^2+1}{V^3},
    \end{aligned}
   \end{equation}
   which is again an element of $C_g$.

Note that $V$ has a strict  minimum $V(0)=1$ and $U$ has a strict maximum $U(0)=n$.
Thus, we see that $g(r)= \omega_nr^n$, and hence  that $V$ is harmonic.
This is a contradiction by the maximum principle of harmonic functions,
 which completes the proof of Theorem 7.

 \vskip 0.50cm
\noindent{\bf ACKNOWLEDGMENT.} The authors would like to express
their deep thanks  to Professor Young Ju Lee and  Professor Do Yong Kwon for valuable advices that improved
this article.

\vskip 1.0 cm


\vskip0.5cm



\begin{thebibliography}{5.4}

\bibitem {a}
 Axler, S.,  Bourdon, P. and  Ramey, W., {\it Harmonic function theory. 2nd ed.,
  Graduate Texts in Mathematics 137}, Springer-Verlag, New York, 2001.

\bibitem {b}
 Blaschke, W., {\it Vorlesungen \"uber Differentialgeometrie und geometrische Grundlagen von Einsteins Relativitatstheorie,
  Band I, Elementare Differentialgeometrie,
  (German) 3rd ed.}, Dover Publications, New York, N. Y., 1945.


\bibitem {h}
 Hartman, P., {\it On complete hypersurfaces of nonnegative sectional curvatures and constant $m$th mean curvature},
  Trans. Amer. Math. Soc. 245 (1978), 363-374.

\bibitem {j}
J\"orgens, K.,
{\it
\"Uber die Losungen der Differentialgleichung $rt-s^2=1$},
Math. Ann. 127(1954), 130-134.

\bibitem {kk1}
Kim, D.-S.  and  Kim, Y. H., {\it Some characterizations of spheres and elliptic paraboloids}, Linear Algebra Appl.
437 (2012), no. 1, 113-120.

\bibitem {kk2}
Kim, D.-S.  and  Kim, Y. H., {\it New approach to Archimedean theorems}, Submitted.

\bibitem {po}
Pogorelov, A. V., {\it
On the improper convex affine hyperspheres},
 Geom. Dedicata 1(1972), no. 1, 33-46.


\bibitem {pr}
Pressley, A., {\it Elementary differential geometry,}
Undergraduate Mathematics Series, Springer-Verlag London, Ltd., London, 2001.

\bibitem {r}
Rosenberg, H.,  {\it Hypersurfaces of constant curvature in space forms},
 Bull. Sci. Math. 117 (1993), no. 2, 211-239.


\bibitem {s}
Stamm, O.,
{\it Umkehrung eines Satzes von Archimedes \"uber die Kugel, }
Abh. Math. Sem. Univ. Hamburg 17(1951), 112-132.


\bibitem {se}
Stein, S., {\it Archimedes. What did he do besides cry Eureka?},
Mathematical Association of America, Washington, DC, 1999.

\bibitem {t}
Thorpe, J. A., {\it Elementary topics in differential geometry},
Undergraduate Texts in Mathematics,
Springer-Verlag, New York-Heidelberg, 1979.

\end{thebibliography}
\end{document}